\documentclass[12pt]{article}
\usepackage[intlimits]{amsmath}

\usepackage{amsfonts,amssymb}
\usepackage[cp1251]{inputenc}
\usepackage[dvips]{graphicx}
\DeclareMathOperator{\mo}{mod}
\usepackage{slashbox,multirow,makecell,longtable}

\newtheorem{theorem}{Theorem\lefteqn{\textstyle \bf \,\,\,\,\,.}}
\newtheorem{conj}{Conjecture\lefteqn{\textstyle \bf \,\,\,\,\,.}}
\title{\bf On Erd\H{o}s--Szekeres problem and related problems\footnote{The work is done with the financial support of the grant RFBR 09-01-00294.}}
\author{V.A. Koshelev\\\small e-mail: {\it koshelev@mccme.ru}}
\date{}
\begin{document}

\voffset=-20mm
\hoffset=-12mm
\font\Got=eufm10 scaled\magstep2 \font\Got=eufm10

\renewcommand{\refname}{References}

\renewcommand{\figurename}{Picture}

\renewcommand{\tablename}{Table}

\renewcommand{\abstract}{\small {\bf Abstract:}}

\newcommand{\pic}[3]{
\begin{figure}[!ht]
\begin{center}
\includegraphics[scale=#3]{#1.eps}
\end{center}
\caption{\footnotesize #2} \label{#1}
\end{figure}
}

\maketitle

\begin{abstract}
Here we give a short survey of our new results. References to the complete
proofs can be found in the text of this article and in the litterature.
\end{abstract}

\section{Introduction and statement of problems}

In 1935 Paul Erd\H{o}s and George Szekeres formulated the following problem  (see \cite{ES}, \cite{Low}).
\vskip+0.2cm
{\bf First Erd\H{o}s--Szekeres problem.} {\it For any integer $n\ge 3$, find the minimal positive number $g(n)$ such that any planar set of points in general position containing at least $g(n)$ points has a subset of cardinality $n$ whose elements are the vertices of a convex $n$-gon.}
\vskip+0.2cm
In 1978 Erd\H{o}s suggested the following modification of the first problem (see \cite{E}).
\vskip+0.2cm
{\bf Second Erd\H{o}s--Szekeres problem.} {\it For any integer $n\ge 3$, find a minimal positive number $h(n)$ such that any planar set ${\cal X}$ in general position containing at least $h(n)$ points has a subset of cardinality $n$ whose elements are the vertices of an empty convex $n$-gon, i.e., of an $n$-gon containing no other points of ${\cal X}$.}
\vskip+0.2cm

Recall that a set of points on the plane is in the {\it general position} if any three of its elements do not lie in a straight line.

The above problems are classical in combinatorial geometry and Ramsey theory (see \cite{TRam}, \cite{Hall}, \cite{Sol}, \cite{Ram}). They can both be generalized as follows.
\vskip+0.2cm

{\bf Third Erd\H{o}s--Szekeres-type problem.} {\it For any integers $n\ge 3$ and $k\ge 0$, find a minimal positive number $h(n,k)$ such that any planar set ${\cal X}$ in general position containing at least $h(n,k)$ points has a subset of cardinality $n$ whose elements are the vertices of convex $n$-gon $C$ with $|(C\setminus \partial C)\cap {\cal X}|\le k$; i.e., the interior of this $n$-gon contains at most $k$ other points of ${\cal X}$.}
\vskip+0.2cm

One more generalization was suggested in \cite{BDV} by Bialostocki, Dierker, and Voxman.

{\bf Fourth Erd\H{o}s--Szekeres-type problem.} {\it For any integers $n\ge 3$ and $q\ge 0$, find a minimal positive number $h(n,\mo q)$ such that any planar set ${\cal X}$ in general position containing at least $h(n,\mo q)$ elements has a subset of cardinality $n$ whose elements are the vertices of convex $n$-gon $C$ with $|(C\setminus \partial C)\cap {\cal X}| \equiv 0 \pmod q$; i.e., the interior of this $n$-gon contains other points from ${\cal X}$ and their number is a multiple of $q$.}

One may find more detailed history of Erd\H{o}s -- Szekeres problems, for example, in the following surveys \cite{BK}, \cite{BMP}, \cite{Sol}.


\section{On the first and second problems}

The first problem was considered by Erd\H{o}s and Szekeres in the article \cite{ES}.
They proved the existence of $g(n)$ for arbitrary $n$ by demonstrating the upper bound $g(n)\le{2n-4\choose n-2}+1$, and they gave the following conjecture:

\begin{conj}
$g(n)=2^{n-2}+1$.
\label{Low_conj}
\end{conj}

This conjecture is proved for $n\le6$. The case $g(3)=3$
is obvious here; equality $g(4)=5$ was proved by E. Klein in 1935
(see pic. \ref{E_Klein},where all three essentially different ways of placing five points on the plane are displayed);
expression $g(5)=9$ was obtained by E. Makai (see \cite{ES},
\cite{Low}, \cite{Sol}); the fact $g(6)=17$ was established rather recently by G. Szekeres, B. McKay and L. Peters in \cite{SL}. Besides, in 1961 Erd\H{o}s and Szekeres have also proved the lower bound $g(n)\ge 2^{n-2}+1$ (see
\cite{Low}).

\pic{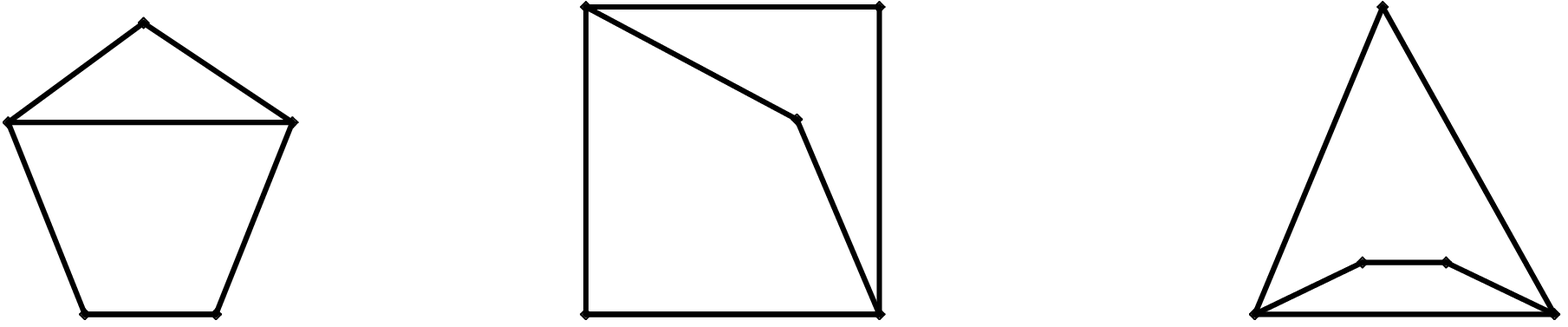}{Any set of five points contains a convex quadrilateral}{0.7}

Inequality $g(n)\le{2n-4\choose n-2}+1$ was repeatedly improved.
The strongest result was obtained in 2005 by G. Toth and P. Valtr: $g(n)\le{2n-5\choose n-3}+1$ (here $ n \ge 5 $;
see \cite {TV05}). Thus, the Erd\H{o}s -- Szekeres conjecture is still neither proved nor disproved, and it is only known that $$2^{n-2}+1 \le g(n)\le{2n-5\choose n-3}+1.$$

In connection with bounding $g(n)$ Erd\H{o}s and Szekeres introduced the notions of {\it cup} and {\it cap}. We assume that a coordinate system $(x, y)$ is fixed in the plane. Let ${\cal X}=\{(x_1,y_1), (x_2,y_2),\ldots, (x_m,y_m) \}$ be a set of points in general position in the plane, with $x_i\neq x_j$ for all $i\neq j$. A subset of points $\{(x_{i_1},y_{i_1}),(x_{i_2},y_{i_2}),\ldots,(x_{i_r},y_{i_r})\}$ is called an {\it $r$-cup} (see pic. \ref{cupcap}) if $x_{i_1}<x_{i_2}<\ldots<x_{i_r}$ and
$$\frac{y_{i_1}-y_{i_2}}{x_{i_1}-x_{i_2}}<\frac{y_{i_2}-y_{i_3}}{x_{i_2}-x_{i_3}}<\ldots<\frac{y_{i_{r-1}}-y_{i_r}}{x_{i_{r-1}}-x_{i_r}};$$

\pic{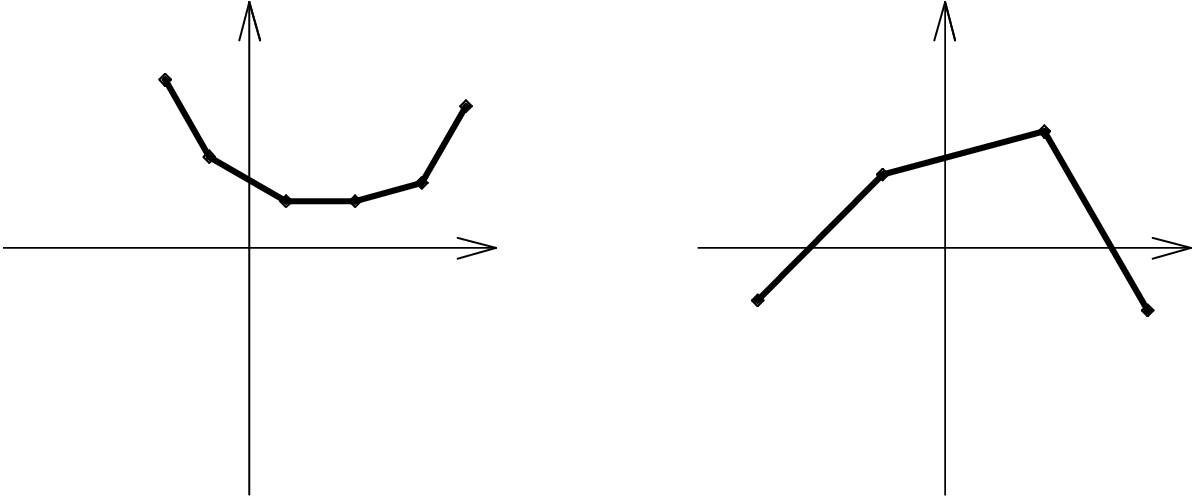}{cup and cap}{0.65}

Similarly, the subset is called an {\it $r$-cap} (see pic. \ref{cupcap}) if $x_{i_1}<x_{i_2}<\ldots<x_{i_r}$ and
$$\frac{y_{i_1}-y_{i_2}}{x_{i_1}-x_{i_2}}>\frac{y_{i_2}-y_{i_3}}{x_{i_2}-x_{i_3}}>\ldots>\frac{y_{i_{r-1}}-y_{i_r}}{x_{i_{r-1}}-x_{i_r}}.$$

Define $f(l,m)$ to be the smallest positive integer for which ${\cal X}$ contains an $l$-cup
or an $m$-cap whenever ${\cal X}$ has at least $f(l,m)$ points.

The problem of finding $f(l,m)$ was completely solved by Erd\H{o}s and Szekeres (see \cite{ES},\cite{Low}). They proved that $f(l,m)={{l+m-4}\choose{l-2}}+1$. Note that the first bound of $g(n)$ is based on the inequality $g(n)\le f(n,n)={{2n-4}\choose{n-2}}+1$. 

The second problem is more deeply understood. Thus, equalities $h(3)=3$
and $h(4)=5$ for it are obvious (see pic. \ref{E_Klein}).
Expression $h(5)=10$ was obtained by H. Harborth in 1978 (see \cite{Harb}). And in 1983 J. Horton proved that $h(n)$ does not exist where $n\ge7$ (see \cite{Hort}). Actually, Horton proved non-existence of $h(n,0)$ where $n\ge 7$. The question of existence and exact value of $h(6)$ (or, which is the same, $h(6,0)$) has been remaining open for a long time. Only in 2006 T. Gerken proved the existence of $h(6)$, by demonstrating the upper bound $h(6)\le g(9)\le {13\choose 6}+1=1717$
(see \cite{Gerken}). Independently of him, C. Nicolas (see \cite{Nic}) and Valtr (see \cite{Val}) presented their proofs, but their upper bounds are worse and equal to, respectively, $g(25)$ and $g(15)$. In 2007 the upper bound was improved by the author of this article:

\begin{theorem}
$h(6)\le 463$ (see \cite{Doclady1},\cite{ENDM},\cite{MSb}).
\end{theorem}

The trivial lower bound $h(6)\ge g(6)\ge 17$ is a consequence of one of the Erd\H{o}s--Szekeres theorems (see \cite{Low}). All the other lower bounds for $h(6)$ were obtained by the computer search. The first one of them was given by D.Rappaport in 1985: $h(6)\ge 21$ (see \cite{Rap}); the second one was done by M.Overmars, B. Scholten and I. Vincent in 1988: $h(6)\ge27$ (see.
\cite{Ov88}). The best known lower bound was obtained in 2003 by Overmars: $h(6)\ge30$ (see \cite{Over}).
Thus, for $h(6)$ estimates $30\le h(6)\le 463$ are proved at present.

\section{On the third problem}

As it is easy to see, for the third problem inequalities $g(n)$ $\le h(n,k)\le h(n)$ are always correct if the appropriate expressions exist. Moreover, $h(n)=h(n,0)\ge h(n,1)\ge h(n,2)\ge h(n,3)\ge \ldots$ and there is a $k'$ such that $h(n,k)=g(n)$ for all $k\ge k'$. For small values of $n$ the following results are obvious: $h(3,k)=3$, $h(4,k)=5$, $h(5,0)=10$, $h(5,\ge 1)=9$. The last result follows from the fact that a convex pentagon with two or more points inside always contains a smaller convex pentagon.

Some results relating to the third problem are obtained in an article by Bl. Sendov (see \cite{Sen}). In this article, with the use of the Horton set (see \cite{Hort}), through which the non-existence of  $h(7)$ was proved, non-existence of $h(n,k)$ was proved for certain values of $k$ where $n>7$. More precisely, $k$ should be less than or equal to $(r+4)2^{m-1}-4m-r-1$, provided $n+2=4m+r$, where $m$ is integer and $r\in\{0,1,2,3\}$. The similar results are obtained in the article by H. Nyklova (see \cite{Nyk}), besides it is proved there that  $h(6,\ge 6)=g(6)$ and the result $h(6,5)=19$ is presented. Note that Sendov's and Nyklova's estimates are asymptotically equal to $(\sqrt[4]{2}+o(1))^n$.

With respect to the fact that all results for $g(6)$ and $h(6)$ were obtained rather recently, the study of the value $h(6,1)$ is interesting (values of $k$, other than 1 may be not so interesting with respect to the conjecture set forth below). We found the upper bound for $h(6,1)$ much better than the upper bound for $h(6,0)$.

\begin{theorem}
The inequality holds $h(6,1)\le g(7) \le 127$ (see \cite{FPM}).
\label{th_h_6_1}
\end{theorem}

Thus, it appears that at present the estimates $17\le h(6,1)\le 127$ are proved. Note that, if the conjecture \ref{Low_conj} of Erd\H{o}s and Szekeres is true, the equality in Theorem \ref{th_h_6_1} will look as $h(6,1)\le g(7)=33$.

Actually, we suppose that the stronger statement is true:

\begin{conj}
$h(6,1)=g(6)=17$.
\end{conj}

Note that it follows immediately from the conjecture that $h(6,1)=h(6,2)=h(6,3)=h(6,4)=h(6,5)=17$. The supposed equality $h(6,5)=17$ obviously contradicts the result of Nyklova set forth above. The point is that this result was proved inaccurately and there are counterexamples to it.

Now we formulate a new result on the existence of $h(n,k)$ for all $n$.

\begin{theorem}
For odd and for even $n$ respectively, the following values do not exist\\ $h(n,{{n-7}\choose{(n-7)/2}}-1)$, $h(n,2{{n-8}\choose{(n-8)/2}}-1)$ (see \cite{ENDM2009},\cite{MZm2}).
\end{theorem}

Note that this theorem gives an asymptotic lower estimate of the form $(2+o(1))^n$ for the maximal value of $k$ such that $h(n,k)$ does not exist. This result is much better than the above-mentionted result of Sendov (see \cite{Sen}). In table \ref{tab_h_nexist} we compare maximum values of $k$ such that $h(n,k)$ does not exist according to Sendov (see \cite{Sen}), Nyklova (see \cite{Nyk}) and this author.

\begin{table}[ht]
\centering
\tabcolsep=2pt
\begin{small}
\begin{tabular}{|r||*{19}{c|}}
\hline
$n, k \ge $ & 7 & 8 & 9 & 10 & 11 & 12 & 13 & 14 & 15 & 16 & 17 & 18 & 19 & 20 & 21 & 22 & 23 & 24 & 25 \\[2mm]\hline
Bl. Sendov 1995 & 0 & 1 & 2 & 3 & 6 & 9 & 12 & 15 & 22 & 29 & 36 & 43 & 58 & 73 & 88 & 103 & 134 & 165 & 196 \\[2mm]\hline
H. Nyklova 2000 & 0 & 1 & 2 & 3 & 6 & 9 & 13 & 19 & 27 & 39 & 51 & 63 & 91 & 119 & 147 & 175 & 238 & 301 & 373 \\[2mm]\hline
V. Koshelev 2009 & 0 & 1 & 2 & 3 & 6 & 11 & 19 & 39 & 69 & 139 & 251 & 503 & 923 & 1847 & 3431 & 6863 & 12869 & 25739 & 48619 \\[2mm] \hline
\end{tabular}
\end{small}
\caption{Comparing lower bounds for $k$ \label{tab_h_nexist}}
\end{table}

It is of interest to find such values of $k$ that $h(n,k)=g(n)$ or $h(n,k)>g(n)$. However we do not know the exact values of $g(n)$. We only know the conjecture \ref{Low_conj}. So we will prove an estimate concerning the maximum value of $k$ for which $h(n,k)>2^{n-2}+1$.

\begin{theorem}
If $n\ge 6$, then $h\left(n,{(n-3)\choose {\lceil(n-3)/2\rceil}}-\left\lceil{\frac{n}{2}}\right\rceil\right)>2^{n-2}+1$ (see \cite{MZm2}).
\end{theorem}

\begin{table}[ht]
\centering
\tabcolsep=1pt
\begin{small}
\begin{tabular}{|r||*{20}{c|}}
\hline
$n$ & 6 & 7 & 8 & 9 & 10 & 11 & 12 & 13 & 14 & 15 & 16 & 17 & 18 & 19 & 20 & 21 & 22 & 23 & 24 & 25 \\[2mm]\hline
$k$ & 0 & 2 & 6 & 15 & 30 & 64 & 120 & 245 & 455 & 916 & 1708 & 3423 & 6426 & 12860 & 24300 &
48609 & 92367 & 184744 & 352704 & 705419 \\ \hline
\end{tabular}
\end{small}
\caption{Values of $k$ such that $h(n,k)>2^{n-2}+1$ \label{tab_h_neq}}
\end{table}

We give some similar results for cups and caps. We define $f(l,m,l_1,m_1)$ to be the smallest positive integer such that any set ${\cal X}$ in general position with no two points having the same x-coordinate and of cardinality $f(l,m,l_1,m_1)$ contains a $l$-cup with at most $l_1$ points inside or an $m$-cap with at most $m_1$ points inside.

\begin{table}[!ht]
\centering
\tabcolsep=3pt
\begin{small}
\begin{tabular}{|c||*{10}{c|}}
\hline
\backslashbox{l}{m} & 5 & 6 & 7 & 8 & 9 & 10 & 11 & 12 & 13 \\[2mm]\hline\hline

5 & \makecell[c]{0 0 \\}& \makecell[c]{0 1 \\}& \makecell[c]{0 4 \\}& \makecell[c]{0 7 \\}& \makecell[c]{0 14 \\}& \makecell[c]{0 21 \\}& \makecell[c]{0 36 \\}& \makecell[c]{0 51 \\}& \makecell[c]{0 82 \\}\\ \hline
6 & \makecell[c]{1 0 \\}& \makecell[c]{1 1 \\}& \makecell[c]{1 4 \\2 2 \\}& \makecell[c]{1 9 \\2 5 \\3 3 \\}& \makecell[c]{1 15 \\2 14 \\3 7 \\4 4 \\}& \makecell[c]{1 29 \\2 23 \\3 19 \\4 9 \\5 5 \\}& \makecell[c]{2 44 \\3 31 \\4 24 \\5 11 \\6 6 \\}& \makecell[c]{1 73 \\2 65 \\3 59 \\4 39 \\5 29 \\6 13 \\7 7 \\}& \makecell[c]{2 110 \\3 87 \\4 74 \\5 47 \\6 34 \\7 15 \\8 8 \\}\\ \hline
7 & \makecell[c]{4 0 \\}& \makecell[c]{4 1 \\2 2 \\}& \makecell[c]{5 5 \\}& \makecell[c]{2 14 \\5 11 \\9 9 \\}& \makecell[c]{5 29 \\9 19 \\14 14 \\}& \makecell[c]{9 49 \\14 29 \\20 20 \\}& \makecell[c]{5 89 \\9 79 \\14 74 \\20 41 \\27 27 \\}& \makecell[c]{9 149 \\14 119 \\20 104 \\27 55 \\35 35 \\}& \makecell[c]{14 224 \\20 167 \\27 139 \\35 71 \\44 44 \\}\\ \hline
8 & \makecell[c]{7 0 \\}& \makecell[c]{9 1 \\5 2 \\3 3 \\}& \makecell[c]{14 2 \\11 5 \\9 9 \\}& \makecell[c]{19 19 \\}& \makecell[c]{9 49 \\19 39 \\34 34 \\}& \makecell[c]{19 99 \\34 69 \\55 55 \\}& \makecell[c]{34 174 \\55 111 \\83 83 \\}& \makecell[c]{19 299 \\55 279 \\83 167 \\119 119 \\}& \makecell[c]{34 524 \\55 447 \\83 419 \\119 239 \\164 164 \\}\\ \hline
9 & \makecell[c]{14 0 \\}& \makecell[c]{15 1 \\14 2 \\7 3 \\4 4 \\}& \makecell[c]{29 5 \\19 9 \\14 14 \\}& \makecell[c]{49 9 \\39 19 \\34 34 \\}& \makecell[c]{74 14 \\14 74 \\69 69 \\}& \makecell[c]{34 174 \\69 139 \\125 125 \\}& \makecell[c]{69 349 \\125 251 \\209 209 \\}& \makecell[c]{125 629 \\209 419 \\329 329 \\}& \makecell[c]{209 1049 \\329 659 \\494 494 \\}\\ \hline
10 & \makecell[c]{21 0 \\}& \makecell[c]{29 1 \\23 2 \\19 3 \\9 4 \\5 5 \\}& \makecell[c]{49 9 \\29 14 \\20 20 \\}& \makecell[c]{99 19 \\69 34 \\55 55 \\}& \makecell[c]{174 34 \\139 69 \\125 125 \\}& \makecell[c]{279 55 \\55 279 \\251 251 \\}& \makecell[c]{125 629 \\251 503 \\461 461 \\}& \makecell[c]{251 1259 \\461 923 \\791 791 \\}& \makecell[c]{461 2309 \\791 1583 \\1286 1286 \\}\\ \hline
11 & \makecell[c]{36 0 \\}& \makecell[c]{44 2 \\31 3 \\24 4 \\11 5 \\6 6 \\}& \makecell[c]{89 5 \\79 9 \\74 14 \\41 20 \\27 27 \\}& \makecell[c]{174 34 \\111 55 \\83 83 \\}& \makecell[c]{349 69 \\251 125 \\209 209 \\}& \makecell[c]{629 125 \\503 251 \\461 461 \\}& \makecell[c]{1049 209 \\209 1049 \\923 923 \\}& \makecell[c]{461 2309 \\923 1847 \\1715 1715 \\}& \makecell[c]{923 4619 \\1715 3431 \\3002 3002 \\}\\ \hline
12 & \makecell[c]{51 0 \\}& \makecell[c]{73 1 \\65 2 \\59 3 \\39 4 \\29 5 \\13 6 \\7 7 \\}& \makecell[c]{149 9 \\119 14 \\104 20 \\55 27 \\35 35 \\}& \makecell[c]{299 19 \\279 55 \\167 83 \\119 119 \\}& \makecell[c]{629 125 \\419 209 \\329 329 \\}& \makecell[c]{1259 251 \\923 461 \\791 791 \\}& \makecell[c]{2309 461 \\1847 923 \\1715 1715 \\}& \makecell[c]{3959 791 \\791 3959 \\3431 3431 \\}& \makecell[c]{1715 8579 \\3431 6863 \\6434 6434 \\}\\ \hline
13 & \makecell[c]{82 0 \\}& \makecell[c]{110 2 \\87 3 \\74 4 \\47 5 \\34 6 \\15 7 \\8 8 \\}& \makecell[c]{224 14 \\167 20 \\139 27 \\71 35 \\44 44 \\}& \makecell[c]{524 34 \\447 55 \\419 83 \\239 119 \\164 164 \\}& \makecell[c]{1049 209 \\659 329 \\494 494 \\}& \makecell[c]{2309 461 \\1583 791 \\1286 1286 \\}& \makecell[c]{4619 923 \\3431 1715 \\3002 3002 \\}& \makecell[c]{8579 1715 \\6863 3431 \\6434 6434 \\}& \makecell[c]{15014 3002 \\3002 15014 \\12869 12869 \\}\\ \hline

\end{tabular}
\end{small}
\caption{Values of $l_1 $ and $m_1$ such that $f(l,m,l_1,m_1)$ does not exist \label{tab_f_nexist}}
\end{table}

\begin{theorem}
Let $c(r)=2^{\lfloor\frac{r-2}{2}\rfloor}+2^{\lceil\frac{r-2}{2}\rceil}-r-1$. If for
$l_0$ and $m_0$ we have $c(l_0)>0,c(m_0)>0$, then
for any $l\ge 5$ and $m\ge 5$, with $l\ge l_0$ and $m\ge m_0$, the following value does not exist:
$f\left(l,m,c(l_0){{l+m-l_0-m_0}\choose {l-l_0}}-1,c(m_0){{l+m-l_0-m_0}\choose {m-m_0}}-1\right)$ (see \cite{MZm2}).
\label{th_f_nexist}
\end{theorem}

\begin{theorem}
For any $l\ge 4$, $m\ge 4$ and $ a \ge 0 $ (with all non-negative arguments of $f$) the inequalities hold
$f\left(l,m,{{l+m-6}\choose {l-3}}-m+1,a\right)>f(l,m)$,
$f\left(l,m,a,{{l+m-6}\choose {m-3}}-l+1\right)>f(l,m)$ (see \cite{MZm2}).
\label{th_f_neq}
\end{theorem}

Tables \ref{tab_f_nexist} and \ref{tab_f_neq} consist of numbers illustrating Theorem \ref{th_f_nexist} and Theorem \ref{th_f_neq}. 

\begin{table}[ht]
\centering
\tabcolsep=4pt
\begin{small}
\begin{tabular}{|c||*{13}{c|}}
\hline
\backslashbox{l}{m} & 4 & 5 & 6 & 7 & 8 & 9 & 10 & 11 & 12 & 13 & 14 & 15 \\[2mm]\hline\hline
4 & --- & --- & --- & --- & --- & --- & --- & --- & --- & --- & --- & --- \\ \hline
5 & 0 & 2 & 5 & 9 & 14 & 20 & 27 & 35 & 44 & 54 & 65 & 77 \\ \hline
6 & 1 & 6 & 15 & 29 & 49 & 76 & 111 & 155 & 209 & 274 & 351 & 441 \\ \hline
7 & 2 & 11 & 30 & 64 & 119 & 202 & 321 & 485 & 704 & 989 & 1352 & 1806 \\ \hline
8 & 3 & 17 & 51 & 120 & 245 & 454 & 783 & 1277 & 1991 & 2991 & 4355 & 6174 \\ \hline
9 & 4 & 24 & 79 & 204 & 455 & 916 & 1707 & 2993 & 4994 & 7996 & 12363 & 18550 \\ \hline
10 & 5 & 32 & 115 & 324 & 785 & 1708 & 3423 & 6425 & 11429 & 19436 & 31811 & 50374 \\ \hline
11 & 6 & 41 & 160 & 489 & 1280 & 2995 & 6426 & 12860 & 24299 & 43746 & 75569 & 125956 \\ \hline
12 & 7 & 51 & 215 & 709 & 1995 & 4997 & 11431 & 24300 & 48609 & 92366 & 167947 & 293916 \\ \hline
13 & 8 & 62 & 281 & 995 & 2996 & 8000 & 19439 & 43748 & 92367 & 184744 & 352703 & 646632 \\ \hline
14 & 9 & 74 & 359 & 1359 & 4361 & 12368 & 31815 & 75572 & 167949 & 352704 & 705419 & 1352064 \\ \hline
15 & 10 & 87 & 450 & 1814 & 6181 & 18556 & 50379 & 125960 & 293919 & 646634 & 1352065 & 2704142 \\ \hline
\end{tabular}
\end{small}
\caption{Values of $l_1$ such that $f(l,m,l_1,?)>f(l,m)$ \label{tab_f_neq}}
\end{table}

\section{On the fourth problem}

~~~Concerning the fourth problem, Bialostocki, Dierker, and Voxman conjectured that $h(n,\mo q)$ exists for all $n\ge 3$ and $q\ge 2$. This conjecture has neither been proved nor disproved thus far. Below, we present the results available and their improvements.

For small values of $n$ the following results are obvious: $h(3,\mo q)=3$, $h(4,\mo q)=5$, $h(5,\mo q)=h(5)=10$. If $n=6$, then $h(6,\mo q)=h(6)$ for all $q$, except finite set of values. Probably $h(6,\mo 2)=g(6)=17$, but for $q\ge 3$ it is possible to construct a point set that shows that $h(n,\mo q)>17$.

Bialostocki, Dierker, and Voxman proved their conjecture (see \cite{BDV}) for $n\ge q+2$ and obtained the upper
bound $h(n,\mo q) \le g\left(R_3(\,\underbrace{n',n',\ldots,n'}_q\,)\right)$, where $n'$ is the minimum positive integer satisfying $n'\ge n$ and $n'\equiv 2 (\mo q)$. Here, $ R_k(l_1, \ldots, l_s) $ is the Ramsey number for complete $k$-uniform hypergraphs with edges painted in s colors in which at least one monochromatic $l_i$-clique with suitable $i$ is sought (see \cite{TRam}, \cite{Hall}). In last formula, the Ramsey number has $q$ arguments with a value of $n'$. Only astronomical estimates of Ramsey numbers are known. In this case, we have a tower of exponentials.

In 1996, Caro (see \cite{Caro}) obtained a more general result for points in the plane with assigned values from a finite Abelian group and for convex polygons with a zero inside sum. As applied to the problem under discussion, his theorem gives $h(n,\mo q)\le 2^{c(q)n}$. Here, $c(q)$ is a function independent of $n$ but growing superexponentially in $q$. Thus, we again deal with a multiple exponential.

Of course, last bound has to be refined. This can be done in two directions. On the one hand, it would be desirable to get rid of the superexponential bounds at least for some relations between $n$ and $q$. On the other hand, the constraint $n\ge q+2$, under which $h(n,\mo q)$ always exists, seems excessive.

The only result in the first direction was obtained in \cite{KPT}, namely, $h(n, \mo q)$ exist for $n\ge 5q/6+O(1)$, but the upper bound is even worse than in \cite{BDV}. Note that a similar result with $n\ge 3q/4+O(1)$ was announced by Valtr, but this result was not published. In the second direction, new results have not been obtained at all. Caro conjectured that $ h(n,\mo q) \le g(c(q)+n) $ with some $c(q)$. We managed to prove the following result.

\begin{theorem}
If $n\ge 2q-1$, then $h(n,\mo q)\le g(q(n-4)+4)$ (see \cite{Doclady2},\cite{ENDM2009},\cite{MZm}).
\label{modq}
\end{theorem}

This theorem considerably improves Caro estimate, since $ g(q(n-4)+4) \le 2^{2qn+O(1)} $. Thus, we eventually have got rid of the multiple exponentials in the inequalities.

However, the constraint $ n \ge 2q-1 $ is somewhat stronger than before, and Caro's conjecture has not been proved (or disproved). Nevertheless, this is an important step toward the solution of the problem.

Note that Bialostocki–-Dierker–-Voxman estimate admits a fairly curious refinement, which is weaker than Caro's result and Theorem \ref{modq}, but, in our view, deserves to be mentioned.

\begin{theorem}
$h(n,\mo q)\le R_3(\underbrace{n,n,\ldots,n}_{q})$, for even $q$,\\  $h(n,\mo q)\le R_3(g(n),\underbrace{n,\ldots,n}_{q-1})$, for odd $q$ (see \cite{Doclady2},\cite{ENDM2009},\cite{MZm}).
\end{theorem}

The theorem is easy to prove by modifying the original
Bialostocki–-Dierker–-Voxman argument.

\section{Chromatic variant of problems}

\pic{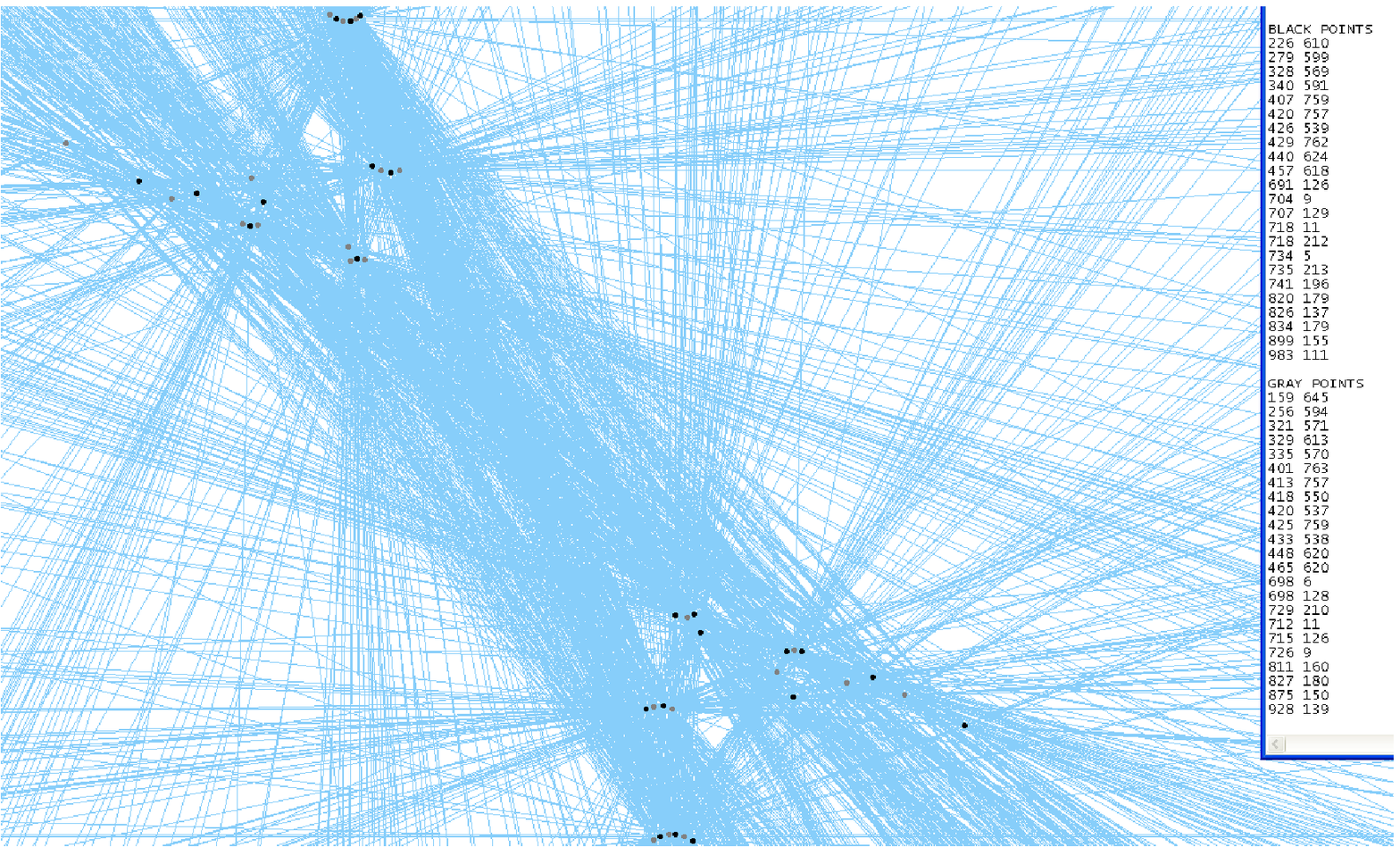}{Example of 46 points}{0.9}

Devillers, Hurtado, K\'arolyi, and Seara \cite{chrom} conjectured that every large enough two-colored set of points, with no three points collinear, contains a convex empty monochromatic fourgon. This can be answered in the affirmative if we omit the condition of convexity (see \cite{4chr}). An example of 18 points with no empty monochromatic convex fourgon from \cite{chrom} led to the problem of finding the maximum number of two-colored points that do not contain an empty monochromatic convex fourgon. Improved lower bounds were given by Brass –- 20 points (see \cite{4chr_20}), Friedman –- 30 points (see \cite{4chr_30}), van Gulik –- 32 points (see \cite{4chr_32}), and finally Huemer and Seara -- 36 points (see \cite{4chr_36}). Here we show (see pic. \ref{46}) a set of 46 two-colored points, no three points collinear, with no empty monochromatic convex fourgons.


\newpage

\begin{thebibliography} {35}

\bibitem{4chr} O. Aichholzer, T. Hackl, C. Huemer, F. Hurtado, and B. Vogtenhuber,
{\it Large bi-colored point sets admit empty monochromatic 4-gons}, 25th European Workshop on Computational Geometry
EuroCG ’09, 25:133–136, Brussels, Belgium, 2009.

\bibitem{BK} I. B\'ar\'any, G. K\'arolyi, {\it Problems and Results around the Erd\H{o}s-Szekeres Convex Polygon Theorem},  Discrete Comput. Geom., LNCS 2089, 91–105, 2001.

\bibitem{BDV} A. Bialostocki, P. Dierker, B. Voxman, {\it Some notes on the Erd\H{o}s--Szekeres theorem},
Discrete Math, 91, 3, 231–238, 1991.

\bibitem{4chr_20} P. Brass, {\it Empty monochromatic fourgons in two-colored point
sets}, Geombinatorics XIV(2), 5–7, 2004.

\bibitem{BMP} P. Brass, W. Moser, J. Pach, {\it Research Problems in Discrete Geometry}, Chapter 8 (Springer, New York, 2005; MCCME, Moscow, 2009).

\bibitem{Caro}  Y. Caro, {\it On the generalized Erd\H{o}s--Szekeres conjecture -- a new upper bound},
Discrete Math, 160, 229–233, 1996.

\bibitem{chrom} O. Devillers, F. Hurtado, G. K\'arolyi, C. Seara, {\it Chromatic
variants of the Erd\H{o}s--Szekeres Theorem}, Computational Geometry,
Theory and Applications 26(3), 193–208, 2003.

\bibitem{ES} P. Erd\H{o}s, G. Szekeres, {\it A combinatorial problem in
geometry}, Compositio Math., 2, 463–470, 1935.

\bibitem{Low} P. Erd\H{o}s, G. Szekeres, {\it On some extremum problems
in elementary geometry}, Ann. Univ. Sci. Budapest E\"otv\"os Sect.
Math., 3–4, 53–62, 1961.

\bibitem{E} P. Erd\H{o}s, {\it Some more problems in elementary geometry},
Austral. Math. Soc. Gaz., 5, 52–54, 1978.

\bibitem{4chr_30} E. Friedman, {\it 30 two-colored points with no empty monochromatic
convex fourgons}, Geombinatorics XIV(2), 53–54, 2004.

\bibitem{Gerken} T. Gerken, {\it On empty convex hexagons in planar
point set}, Discrete Comput. Geom., 39, 239–272, 2008.

\bibitem{TRam} R.L. Graham, B.L. Rothschild, J.H. Spencer,
{\it Ramsey theory}, John Wiley and Sons, NY, Second Edition,
1990.

\bibitem{Hall} M. Hall, Jr., {\it Combinatorial Theory} (Blaisdell, Waltham,
Mass. 1967; Mir, Moscow, 1970).

\bibitem{Harb} H. Harborth, {\it Konvexe F\"unfecke in ebenen Punktmengen},
Elem. Math., 33, 116–118, 1978.

\bibitem{Hort} J.D. Horton, {\it Sets with no empty 7-gons}, Canad. Math.
Bull., 26, 482–484, 1983.

\bibitem{4chr_36} C. Huemer, C. Seara, {\it 36 two-colored points with no empty monochromatic convex fourgons}, Geombinatorics XIX(1), 5–6, 2009.
    
\bibitem{KPT} G. K\'arolyi, J. Pach, G. T\'oth, {\it A modular version of the Erdos-Szekeres theorem},
Studia Sci. Math. Hungar, 38, 245–259, 2001.

\bibitem{Doclady1} V.A. Koshelev, {\it On the Erd\H{o}s--Szekeres problem}, Dokl. Math. 76, 603–605, 2007 [Dokl.
Akad. Nauk 415, 734–736, 2007].

\bibitem{ENDM} V.A.Koshelev, {\it On the Erd\H{o}s--Szekeres problem in combinatorial geometry},
    Electronic Notes in Discrete Mathematics, V. 29, 175–177, 2007.

\bibitem{FPM} V.A. Koshelev, {\it Almost empty hexagons}, Journal of Mathematical Sciences, In press [Fundam. Prikl. Mat. 14, 91–120, 2008].

\bibitem{Doclady2} V.A. Koshelev, {\it Around Erd\H{o}s--Szekeres problems}, Dokl. Math. 79, 360–361, 2009 [Dokl.
Akad. Nauk 426, 304–306, 2009].

\bibitem{ENDM2009} V.A.Koshelev, {\it On Erd\H{o}s--Szekeres-type problems},
Electronic Notes in Discrete Mathematics, V. 34, 447–451, 2009.

\bibitem{MSb} V.A. Koshelev, {\it On Erd\H{o}s -- Szekeres problem for empty hexagons in the plane},
Modeling and analysis of information systems, Vol. 16, 2, 22–74, 2009, On Russian.

\bibitem{MZm} V.A. Koshelev, {\it Erd\H{o}s--Szekeres theorems and congruents}, Mathematical Notes, In press [Matematicheskie Zametki, 87, 2010, In press].

\bibitem{MZm2} V.A. Koshelev, {\it Inside points in Erd\H{o}s--Szekeres theorems}, Mathematical Notes, In press [Matematicheskie Zametki, In press].

\bibitem{Sol} W. Morris, V. Soltan, {\it The Erd\H{o}s - Szekeres problem
on points in convex position}, Bulletin (new series) of the Amer.
Math. Soc., 37, 4, 437–458, 2000.

\bibitem{Nic} C. Nicolas, {\it The empty hexagon theorem}, Discrete Comput. Geom., 38, 2, 389–397, 2007.

\bibitem{Nyk} H. Nyklova, {\it Almost empty polygons},
Studia Scientiarum Mathematicarum Hungarica, 40, 3, 269–286, 2003.

\bibitem{Ov88} M. Overmars, B. Scholten, I. Vincent, {\it Sets
without empty convex 6-gons}, Bull. European Assoc. Theor. Comput.
Sci., 37, 160–168, 1989.

\bibitem{Over} M. Overmars, {\it Finding sets of points without
empty convex 6-gons}, Discrete Comput. Geom., 29, 153–158, 2003.

\bibitem{Ram} F.P. Ramsey, {\it On a problem of formal logic},
Proc. London Math. Soc. Ser. 2, 30, 264–286, 1930.

\bibitem{Rap} D. Rappaport, {\it Computing the largest empty convex subset of a set of points}, ACM 0-89791-163-6/85/006/0161, 161–167, 1985.

\bibitem{Sen} Bl. Sendov, {\it Compulsory configurations of points in the plane}, Fundam. Prikl. Mat. 1, 491–516, 1996, On Russian.

\bibitem{SL} G. Szekeres, L. Peters, {\it Computer solution to the 17-point
Erd\H{o}s--Szekeres problem}, ANZIAM J., 48, 151–164, 2006.

\bibitem{TV05} G. T\'oth, P. Valtr, {\it The Erd\H{o}s--Szekeres
theorem: upper bounds and related results}, Combinatorial and
Computational geometry, MSRI Publication 52, 557–568, 2005.

\bibitem{Val} P. Valtr, {\it On the empty hexagons},
Contemporary Mathematics, 453, 433–442, 2008.

\bibitem{4chr_32} R. Van Gulik, {\it 32 two-colored points with no empty monochromatic
convex fourgons}, Geombinatorics XV(1), 32–33, 2005.

\end{thebibliography}
\end{document}